\documentclass{article}
\usepackage{amsmath,amssymb,color,bbm}
\usepackage{graphicx}
\usepackage{cite}
\usepackage{fullpage}
\usepackage{parskip}
\usepackage{mathrsfs}

\newtheorem{thm}{Theorem}[section]
\newtheorem{lem}[thm]{Lemma}

\newtheorem{obs}[thm]{Observation}

\newtheorem{remark}[thm]{Remark}

\newenvironment{proof}{{\bf Proof.}  }{\hfill$\blacksquare$}

\title{The stability of fixed points for a Kuramoto model with Hebbian interactions.}
\author{Jared C. Bronski\\ Yizhang He \\ Xinye Li \\ Yue Liu \\ Danielle Rae Sponseller \\Seth Wolbert \\University of Illinois \\
Department of Mathematics\\ 1409 W Green St. \\ Urbana, IL 61801. }
\date{\today}

\begin{document}

\maketitle

\begin{abstract} 
We consider a variation of the Kuramoto model with dynamic coupling, where the coupling strengths are allowed to evolve in response to the phase difference 
between the oscillators, a model first considered by  Ha, Noh and Park. In particular we study the stability of fixed points for this model.
We demonstrate a somewhat surprising fact: namely that the fixed points of this model, as well as their stability, can be completely 
expressed in terms of the fixed points and stability of the analogous classical Kuramoto problem where the coupling strengths are fixed to a constant (the same for 
all edges).  In particular for the ``all-to-all'' network, where the underlying graph is the complete graph, the problem reduces to the problem of understanding the fixed points and stability of the all-to-all Kuramoto model  with equal edge weights, a problem that has been completely solved. 

\end{abstract}

{\bf The Kuramoto model is a widely studied model for many diverse synchronization and phase-locking phenomenon. There have been many extensions to the model to incorporate adaptive coupling strengths. We consider a Kuramoto type model that includes adaptive couplings of a Hebbian type, as would be natural for networks of neurons. We establish that the fixed points of model, as well as their stability, can be reduced to the fixed points and stability for the classical Kuramoto model with couplings of fixed strength. This is surprising, given that the 
dimension of the dynamical system for the adaptive model ($N(N+1)$ for the all-to-all topology) greatly exceeds that of the classical Kuramoto model ($N$ for the all-to-all topology), and means that the many results available for the classical Kuramoto model are available for the model with adaptive coupling. In the case of a general nonlinear interaction (not the $\sin(\theta_i-\theta_j)$ interaction) we still reduce the dimension of the stability problem, however the resulting stability problem no longer has an obvious interpretation in terms of the classical Kuramoto model.}
\section{Introduction}

Many natural physical or biological systems such as electrical networks, neural networks, etc, can be modeled as coupled oscillators on some graph or network.\cite{Peskin.75,PRK.book,Medvedev.Kopell.01,Kuramoto.book,Kopell.Ermentrout.86,E,Sync.book}  One fundamental question in 
this setting is whether the oscillators synchronize, or phase lock. One important mathematical model for this phenomenon is the Kuramoto model of coupled oscillators\cite{K}, which can 
be expressed as the following set of coupled ordinary differential equations
\begin{align} 
   \frac{d\theta_{i}}{dt} = {}& \omega_{i} + \sum_{j} \gamma_{ij}\sin(\theta_{j} - \theta_{i})
\label{eqn:Kur}
\end{align}
where $\gamma_{ij}$ is the coupling between oscillators $i$ and $j$ and is assumed to be constant, and $\omega_i$ is the natural frequency of the $i^{th}$ oscillator. 
This model has received a great deal of attention in the literature, and is an important simplified model for a whole host of different phenomena. The Kuramoto model does not, 
however, include an important mechanism known as plasticity. Plasticity is 
the term used in neuroscience to describe the fact that the coupling between 
two neurons is observed to change dynamically in response to the states of 
the  neurons. In particular most types of neurons follow what is known as a 
Hebbian type dynamics, where the coupling between cells increases if they fire 
at close to the same time. This was 
memorably expressed by Schatz\cite{Schatz.1992} as ``cells that fire together wire together''. 
Thus this kind of dynamics is a logical extension to the Kuramoto model. 
Several models for this Hebbian type behavior  have been proposed\cite{Vicente.Arenas.Bonilla.1996,MLHBT.2007,Timms.English.2014,Isakov.Mahadevan.2014,Holzel.Krischer.2015,Ha.Noh.Park.2016}, but the one 
which we consider here is the following 
\begin{align} 
\begin{split}
   \frac{d\theta _i}{dt} = {}& \omega _i + \sum_{j}^{} \gamma _{ij} \sin{(\theta _j - \theta _i)}
   \\
   \frac{d\gamma_{ij}}{dt} = {}& \mu\cos(\theta_{i} - \theta_{j})- \alpha\gamma_{ij}.
\end{split}
\label{eqn:Hebb}
\end{align}
This model appears to have been independently proposed in the literature several times. In the recent paper of Ha, Now and Park\cite{Ha.Noh.Park.2016} this is referred to as Model A. The papers of Isakov and Mahadevan\cite{Isakov.Mahadevan.2014} and Timms and English\cite{Timms.English.2014} consider a similar model 
with the additional terms representing stochastic forcing and axonal delay respectively. 
   
The dynamics of $\gamma_{ij}$ in equation (\ref{eqn:Hebb}) reflect Hebbian dynamics: 
the cosine term implies that the coupling between oscillators that are in phase tends to increase, while the coupling between oscillators which are out of phase tends to decrease. Additionally there is a linear damping term to reflect the lossy nature of this process. Note that this system is a gradient flow, and takes the form $\dot \theta_i = -\partial_{\theta_i} H, \dot \gamma_{ij} = -\mu\partial_{\gamma_{ij}} H$, where the Lyapunov function $H$ is given by 
\[
H = -{\bf \theta}\cdot {\bf \omega} + \sum_{ij} \gamma_{ij}\cos(\theta_i-\theta_j) + \frac{\alpha}{2\mu}\sum \gamma_{ij}^2.
\]
We note that this Lyapunov function is essentially the same as the Lyapunov function for the classical Kuramoto model, the only differences being the third term and the fact that the $\gamma_{ij}$ are now dynamical variables: if one looks at the gradient flow in $\bf\theta$, treating $\bf\gamma$ as constant, one gets the classical Kuramoto model. If  one looks at the gradient flow in $\bf\theta$ and $\bf \gamma$ together one gets the Kuramoto model with Hebbian dynamics. 

For the sake of convenience we will briefly summarize the main results of Ha, Noh and Park.
Let  $\Theta = (\theta_1,...,\theta_N)$ and $\Gamma = [\gamma_{ij}]$,
\begin{align*} 
\begin{split}
   D(\Theta(t)) = \max_{1 \le i,j \le N}|\theta_i(t)-\theta_j(t)|
   \\
   D(\dot \Theta(t)) = \max_{1 \le i,j \le N}|\dot \theta_i(t)-\dot \theta_j(t)|
\end{split}
\end{align*}
Then Ha et. al. prove the following results: 
\begin{itemize}
    \item In the case of identical oscillators ($\omega_i=0$),  if the initial 
phase angles  $\theta_i$ are confined to an arc with a length at most $\frac{\pi}{2}$, then for any solution $\Theta = \Theta(t)$, the solution norm decays at least exponentially.
  \item  For nonidentical oscillators strong synchronization does not hold, but the phase diameter $D(\Theta(t))$ is, for large times $t$, of the order of 
the inverse of the initial minimum coupling strength.
\end{itemize}

In what follows we will consider the above model with $\mu=1$, 
\begin{align} 
\begin{split}
   \frac{d\theta_{i}}{dt} = {}& \omega_{i} + \sum_{j\in N(i)} \gamma_{ij}\sin(\theta_{j} - \theta_{i})
   \label{eqn:kmhebb} \\
   \frac{d\gamma_{ij}}{dt} = {}& \cos(\theta_{i} - \theta_{j}) - \alpha \gamma_{ij}
\end{split}
\end{align}
which can be achieved by appropriate rescaling. Here $N(i)$ denotes the set of neighbors of vertex $i$, and we will let $E$ denote the number of edges in the graph.
 We construct all of the fixed points of this model, and show that they can be related 
to the fixed points of the Kuramoto model on the same network with  {\em fixed} constant edge weights. Moreover, and much more surprisingly,  the Morse index and hence the stability of these 
fixed points is determined completely by the Morse index of the corresponding solutions solutions to the Kuramoto model with fixed edge weights.


\section{Main Results}
As is usual with the Kuramoto model we can, by working in the co-rotating frame, assume that $\sum\omega_i=0$. In this case the phase-locked solutions of 
equation (\ref{eqn:kmhebb}) become fixed points, which satisfy the system of equations 
 \begin{align} 
\begin{split}
   & \omega_{i} =-  \sum_{j\in N(i)} \gamma_{ij}\sin(\theta_{j} - \theta_{i})
   \label{eqn:fp} \\
    &\alpha \gamma_{ij}= \cos(\theta_i-\theta_j).
\end{split}
\end{align}
One can consider allowing the angles $\theta_i$ to vary, in which case equations (\ref{eqn:fp})  become a map from ${\mathbb R}^N$ (the angles $\theta_i$) to 
$({\bf \omega},{\bf \gamma})\in {\mathbb R}^{N+E}$.\footnote{More accurately the co-dimension one subspace such that $\sum \omega_i=0$ } The range of this map is clearly compact, since sine and cosine are bounded functions. We can 
algebraically eliminate $\gamma_{ij}$ to give 
\[
 \omega_{i} =- \frac{1}{\alpha} \sum_{j\in N(i)} \cos(\theta_i-\theta_j) \sin(\theta_{j} - \theta_{i}) = - \frac{1}{2\alpha} \sum_{j\in N(i)}  \sin(2(\theta_{j} - \theta_{i})).
\label{eqn:fixedpt}
\]
As ${\bf \theta}$ ranges over all ${\mathbb R}^N$ the range of this map gives the set of feasible frequencies ${\bf \omega}$: those frequencies for which this model has a fixed point. 
Notice that Equation (\ref{eqn:fixedpt}) is exactly the equation for the 
fixed points for the classical Kuramoto model (\ref{eqn:Kur}) on the same 
graph with edge weights of strength $\gamma_{ij}=\frac{1}{2\alpha}$ and 
angles $2\theta_i$. 
This leads to our first observation 
\begin{obs} Consider the extended Kuramoto model with Hebbian dynamics given by equations (\ref{eqn:kmhebb}), and the Kuramoto model on the same graph with fixed edge weights $\frac{1}{2\alpha}$ given by 
\[
\frac{d\theta_i}{dt} = \omega_i + \frac{1}{2\alpha}\sum_{j\in N(i)}\sin(\theta_j-\theta_i).
\]
Then ${\bf \theta}^*$ is a fixed point for the Kuramoto  model (\ref{eqn:Kur}) for frequency ${\bf \omega}$ if and only if $\frac{{\bf \theta}^*}{2}$ is a fixed point of the Hebbian Kuramoto model for frequency ${\bf \omega}$, with coupling strengths $\gamma_{ij}=\frac{1}{\alpha}\cos(\theta_i^*-\theta_j^*)$.

In other words the two models have exactly the same same feasible sets, but the angular spread for the fixed point of the model with Hebbian dynamics is exactly {\em half} of the angular spread for the model with fixed edge weights and no Hebbian dynamics.   
\label{obs}
\end{obs}

Thus the fixed point set of a model with edge weights that are allowed to 
dynamically evolve is exactly the same as the fixed point set of the model 
where the edge weights are fixed to the same constant $\frac{1}{2\alpha}.$ 
Of course the obvious next question is that of dynamic stability for these 
fixed points: the equivalence of the fixed point sets does not imply that
the dynamics or stability properties are the same. In fact we will show that, 
somewhat surprisingly, the stability is the same as the stability for the 
classical Kuramoto model: for $\alpha$ is positive the dimension of the 
unstable manifold to a fixed point of the model with Hebbian interactions 
is the same as the dimension of the unstable manifold to corresponding fixed 
point of the classical Kuramoto model! 

To analyze the stability of the fixed points, we perform a linearization around the fixed points for the system. The generator of the linearized flow is 
the Jacobian $J$, defined by
$$J_{ij}=\frac{\partial f_i}{\partial x_j}$$
where:
\begin{eqnarray*}
f_i  &=&
\left\{
\begin{array}{ll}
       $$ \omega _i + \sum_{j}^{} \gamma _{ij} \sin{(\theta _j - \theta _i)} &  1 \leq i \leq N $$\\ 
        \cos(\theta_{i'} - \theta_{j'})-\alpha  \gamma_{i'j'} &   N+1 \leq i' \leq N+E \\
\end{array} 
\right.  \\
x_i  &=&
\left\{
\begin{array}{ll}
        \theta_i  &   1 \leq i \leq N \\
        \gamma_{i'j'} &  N+1 \leq i' \leq N+E \\
\end{array} 
\right. 
\end{eqnarray*}
Since the Kuramoto model with Hebbian dynamics is a gradient flow it follows immediately that the Jacobian is a symmetric matrix, and that the eigenvalues of the Jacobian are purely 
real. A bit of algebra shows that the matrix J has the following block partitioned form:
\[
\renewcommand*{\arraystretch}{1.25}
J =   \left(\begin{array}{c|c} 
     A & B\\ \hline 
     B^T & C 
     \end{array}  \right)
\]
Where the matrices $A,B,C$ are $N \times N$, $N \times E $ and $E \times E$ 
respectively. Matrix A is the same as the Jacobian matrix for the original Kuramoto model, which takes the form of a graph Laplacian:
\begin{align*} 
A_{ij}=
\left\{
\begin{array}{ll}
       $$ \gamma_{ij}  \cos(\theta_i-\theta_j)  &   i \neq j  \\
       -\sum_{k}^{}  \gamma_{ij}  \cos(\theta_i-\theta_k) &   i=j $$\\
\end{array} 
\right. 
\end{align*}
Note that this form can be further simplified by the fact that, 
at the fixed point, we have  
$\gamma_{ij}= \frac{1}{\alpha} \cos(\theta_i-\theta_j)$, giving 
\begin{align*} 
A_{ij}=
\left\{
\begin{array}{ll}
       \frac{1}{\alpha}  \cos(\theta_i-\theta_j)^2  &   i \neq j  \\
       -\frac{1}{\alpha} \sum_{k} cos(\theta_i-\theta_k)^2 &   i=j \\
\end{array} 
\right. 
\end{align*}
The matrix C has a particularly simple form, namely a scalar multiple of the 
$E\times E$ identity matrix. 
\begin{eqnarray*}
\frac{\partial f_i}{\partial x_j} &=& \frac{\partial(\cos(\theta_{i'}- \theta_{j'})-\alpha  \gamma_{i'j'})}{\partial \gamma_{i'j'}} = -\alpha\delta_{i,i'}\delta_{j,j'}
\end{eqnarray*}
where $\delta_{i,i'}$ is the Kronecker delta. Finally we come to the 
matrix $B$. Again we compute 
\begin{eqnarray*}
\frac{\partial f_i}{\partial x_j} &=& \frac{\partial(\cos(\theta_{i'}- \theta_{j'})-\alpha  \gamma_{i'j'})}{\partial \theta_k} = \sin(\theta_{i'}-\theta_{j'}) (\delta_{j',k} - \delta_{i',k}).
\end{eqnarray*}
Note that $B$ is a standard object from algebraic graph theory: it is the (weighted) incidence matrix for the Laplacian on the graph, with edge weight given $\sin(\theta_{i'}-\theta_{j'}).$

Once we have determined the Jacobian for the flow, we are interested in its 
eigenvalues. In particular we are interested in counting the number of positive and negative eigenvalues, which give us the dimensions of the stable and 
unstable manifolds to the fixed point. At this point we cite a result of 
Haynsworth\cite{Haynsworth} showing that the inertia or index of a matrix is 
additive in the Schur complement. This will be a useful tool for  counting  
the number of positive or negative eigenvalues of the partitioned matrix. 
\begin{lem}[Haynsworth]
Let $M$ be a real symmetric matrix partitioned as follows
\begin{align*}
\renewcommand*{\arraystretch}{1.35}
M = \left(\begin{array}{c|c} A & B \\\hline  B^{T} & C \\ \end{array}\right).
\end{align*}
Assume that the submatrix $C$ is non-singular. Let $n_+(\cdot)$ (resp. $n_-(\cdot)$) denote the number of positive (resp negative) eigenvalues of the matrix. 
\[
n_{\pm}(M) =  n_{\pm}(C) + n_{\pm}(A - BC^{-1}B^{T}).
\]
\end{lem}
Since this result is critical for the current paper and the proof is not difficult we include a short proof. 

\begin{proof}
This is a relatively straightforward calculation using the Schur complement.\footnote{A term originally coined by Haynsworth.} It is easy to check that 
we have the identity for partitioned matrices
\[
\renewcommand*{\arraystretch}{1.35}
\left(\begin{array}{c|c} I_k & -BC^{-1} \\\hline 0 & I_{m-k}\end{array}\right)
\left(\begin{array}{c|c} A & B \\\hline B^T & C \end{array}\right) 
\left(\begin{array}{c|c} I_k & 0 \\\hline -C^{-1} B^T & I_{m-k}\end{array}\right) 
=\left(\begin{array}{c|c} A -BC^{-1}B^T & 0 \\\hline 0 & C \end{array}\right) \]
where $A$ and $C$ are $k\times k$ and $(m-k)\times (m-k)$ real symmetric matrices respectively. By Sylvester's law of inertia\cite{Sylvester} it follows that 
if $M$ is a real symmetric matrix and $U$ is an invertible real matrix then 
\[
n_{\pm}(M) = n_{\pm}(U^T M U) 
\]   
from which the result follows immediately. 
\end{proof}

We can apply the Haynsworth formula to drastically simplify the stability 
calculation. In our case $C$ is $-\alpha I$. The damping coefficient 
$\alpha$ is assumed to be positive, so $n_+(C)=0$, and any potential 
instabilities are due to the second term, $A-BC^{-1}B^T=A + \frac{1}{\alpha}B B^T.$ 
This is easy to 
compute: it is a standard result in algebraic graph theory that, if $B$ is a 
(weighted) incidence matrix, then $B B^T$ is the graph Laplacian whose edge 
weights are the {\em squares} of the weights in the incidence matrix.  
Therefore we have the following expression for $BB^T$
\[
(BB^T)_{ij}=
\left\{
\begin{array}{ll}
       \sum_{k} \sin^2(\theta_k - \theta_i)  &   i = j  \\
       -\sin^2(\theta_i-\theta_j) &   i \neq j \\
\end{array} 
\right. 
\]

Now we can calculate $A-BC^{-1}B^T$:
\begin{eqnarray*}
(A-BC^{-1}B^T)_{ij} &=& 
\frac{1}{\alpha}\left\{
\begin{array}{ll}
       -\sum_{k}^{} \cos(\theta_i-\theta_k)^2 &   i=j \\
    \cos(\theta_i-\theta_j)^2  &   i \neq j  \\
\end{array} 
\right. 
+\frac{1}{\alpha} \left\{
\begin{array}{ll}
       \sum_{k}^{} \sin(\theta_k - \theta_i)^2  &   i = j  \\
       -\sin(\theta_i-\theta_j)^2 &   i \neq j \\
\end{array} 
\right. \\
&=& \frac{1}{\alpha} 
\left\{
\begin{array}{ll}
       - \sum_{k}^{} \cos(2(\theta_i-\theta_k)) &   i=j \\
         \cos(2(\theta_i-\theta_j))  &   i \neq j  \\
\end{array} 
\right. 
\end{eqnarray*}

Combining this with observation (\ref{obs}) we have proven our main result:
\begin{thm}
 Consider Hebbian Kuramoto model given by equations (\ref{eqn:kmhebb}), and 
the Kuramoto model on the same graph with fixed edge weights 
$\frac{1}{2\alpha}$ given by equation (\ref{eqn:Kur}) with $\gamma=\frac{1}{2\alpha}$.  
Then ${\bf \theta}^*$ is a fixed point for the Kuramoto  model (\ref{eqn:Kur}) for frequency ${\bf \omega}$ if and only if $\frac{{\bf \theta}^*}{2}$ is a fixed point of the Hebbian Kuramoto model for frequency ${\bf \omega}$, with coupling strengths $\gamma_{ij}=\frac{1}{\alpha}\cos(\theta_i^*-\theta_j^*)$. Furthermore the dimensions of the unstable manifolds to each of these fixed points is the same. 

\end{thm}

\begin{remark}
We remark that the stability problem for the classical Kuramoto problem has been extensively studied\cite{Bronski.DeVille.Park.2012,C,Bullo.Dorfler.11,HaHaKim,VO1,VO2,Mirollo.Strogatz.05}, and that all of these results can immediately be applied to the Kuramoto model with Hebbian interactions. 

One can do a similar calculation for a general nonlinear coupling term, but
it is only for the $\sin$ interaction (Kuramoto) that one gets such a simple 
result, essentially due to the angle addition formulae. The natural generalization is to do a gradient flow on the Lyapunov function 
\[
H = - {\bf \omega} {\bf \theta} + \sum_{ij} \gamma_{ij} F(\theta_i-\theta_j) - \frac{\alpha}{2\mu} \sum \gamma_{ij}^2
\]
where $F$ is an even, periodic function. The fixed points are given by the solution to 
\begin{align*}
&0=\omega_i - \sum \gamma_{ij}f(\theta_i-\theta_j) &\\
&\gamma_{ij} = \frac{1}{\mu} F(\theta_i-\theta_j) & 
\end{align*}
where $f=F'$ is an odd, periodic function. 
The Jacobian has the same structure, with the off-diagonal blocks given by 
the (weighted) incidence matrix for the graph, and the results of Haynsworth 
reduce the problem to uderstanding positive eigenvalues of the graph Laplacian
\[
J_{ij} = \left\{\begin{array}{c} \kappa_{ij} \qquad i \neq j \\ -\sum_k \kappa_{ik} \qquad i=j \end{array}\right.
\]  
where the edge weights $K_{ij}=f^2(\theta_i-\theta_j) + F(\theta_i \theta_j)f'(\theta_i-\theta_j).$ In the case of the Kuramoto model the fixed point equation 
simplifies due to the angle addition formula for $\sin$ while the Jacobian 
matrix simplifies due to the angle addition formula for $\cos$. Still, the 
Haynesworth theorem reduces the calculation from determining the inertia of an $(N+E)\times(N+E)$ matrix to determining the inertia of an $N\times N$ matrix.
With a dense network such as the all-to-all topology this is a non-trivial 
reduction.   
\end{remark}

\section{Numerics}

In this section we present some numerical experiments on the Hebbian Kuramoto 
system for three oscillators in the all-to-all topology (the graph is $K_3$). 
Since the natural oscillator frequencies $\omega_i$ can, without loss of 
generality, be assumed to sum to $0$ all graphs are drawn in the mean zero 
frequency plane spanned by $x=(1,-1,0)/\sqrt{2}$ and $y = (1,1,-2)/\sqrt{6}$.

All numerics were done with the parameter values $\mu=1$ and $\alpha=.3$. In each case the initial conditions for the angular variables were taken to be $\theta_1(0)=\theta_2(0)=\theta_3(0)=0$ and the initial conditions for the dynamic 
coupling  strengths were taken to be either $\gamma_{12}(0)=\gamma_{13}(0)=\gamma_{23}(0)=1$ or  $\gamma_{12}(0)=\gamma_{13}(0)=\gamma_{23}(0)=3.$ These initial 
conditions were chosen to be the same in order to preserve the $S_3$ 
permutation symmetry of the problem. All numerics were done with Mathematica\cite{Wolfram}: the ordinary differential equations were solved using the 
Mathematica command {\bf NDSolve} and plotted using the appropriate 
graphics command.

\begin{figure}[h]
 \centering
\includegraphics[width = 0.60\columnwidth]{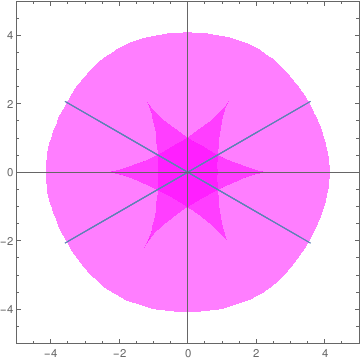} 
\caption{The Region of Stable Phase-Locking for the Hebbian Kuramoto model for $\alpha=.3$}
\label{fig:FeasibleRegion}
\end{figure}

The first figure (Figure (\ref{fig:FeasibleRegion}) depicts the set of 
mean-zero frequencies supporting a 
stable phase-locked solution. In this plot the mean-zero
frequency plane is parameterized by points $(a,b),$ with 
$(\omega_1,\omega_2,\omega_3) = (a/\sqrt{2}+b/\sqrt{6},b/\sqrt{6}-a/\sqrt{2},2b/\sqrt{6}).$  Note that, while the region looks quite circular it actually 
has the symmetry group of the hexagon, $D_3=S_3$. Note that the darker 
star-shaped region in the center of the figure is an artifact of plotting 
with a physical significance. This darker region indicates a multiple covering 
of of the domain and the birth of additional (unstable) solution branches.

The second figure denotes the set of frequencies for which the system actually 
achieved phase-locking, with initial conditions $\theta_{1,2,3}=0$ and two different sets of initial conditions for the coupling strengths: 
$\gamma_{ij}(0)=1$ and $\gamma_{ij}(0)=3$. The numerics in this figure 
were created by numerically integrating the equations out to time $T=75$
and numerically evaluating the $\ell_{\infty}$ norm of the vector field evaluated at the solution point at final time step. 
The four plots in Figure (2-5) are density plots depicting the $\ell_\infty$ norm of the vector field at the final time. Figures (2) and (3) depict the long-time stability for solutions 
to \ref{eqn:Hebb} with two different initial 
values of the coupling constant $\gamma_{ij}(0)=1$ (Figure (2)) and  
$\gamma_{ij}(0)=3$  (Figure (3)).  The dark region depicts the set of 
frequencies for which the size of the vector field  at the terminal time is 
approximately zero, indicating that the solution has converged to a fixed 
point. Note that the set of frequencies for which the solution has converged to a fixed point is noticeably smaller for smaller starting values of the 
coupling strengths $\gamma_{ij}(0).$ Figures (4) and (5) show the same density 
plots overlaid with 
the theoretical result giving the region for which there exists a stable, 
phase-locked solution. We can see for the smaller value of the initial 
coupling constant the solution converges to the fixed point for some 
but not all choices of the initial frequency. We have repeated these numerics 
out to longer times, with no visible change in the result. For larger starting 
values of $\gamma$, on the other hand, the solution seems to converge to the 
fixed point for all values of the frequency for which a fixed point exists. 
Thus it appears that the 
basin of attraction for the fixed point is not the whole phase space. 
This is in contrast to the 
situation for the classical Kuramoto on the complete graph, where the basin 
of attraction for the fixed point appears to be the whole entire phase space 
(excepting a set of measure zero).

\begin{figure}[!htb]
    \centering
    \begin{minipage}{.5\textwidth}
        \centering
        \includegraphics[width=0.85\columnwidth]{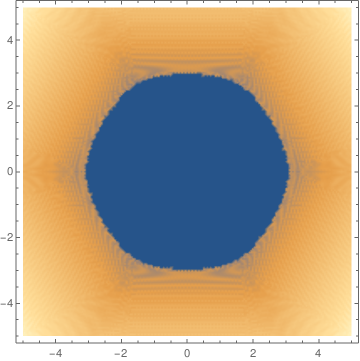}
        \caption{Phase locking for $\gamma_{ij}(0)=1$} 
 \includegraphics[width=0.85\columnwidth]{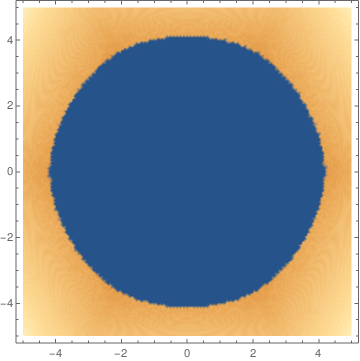}
        \caption{Phase locking for $\gamma_{ij}(0)=3$}
        \label{PhaseLock1}
    \end{minipage}%
    \begin{minipage}{.5\textwidth}
        \centering
        \includegraphics[width=0.85\columnwidth]{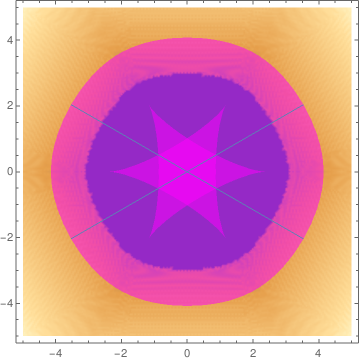}
        \caption{Phase locking for $\gamma_{ij}(0)=1$}
 \includegraphics[width=0.85\columnwidth]{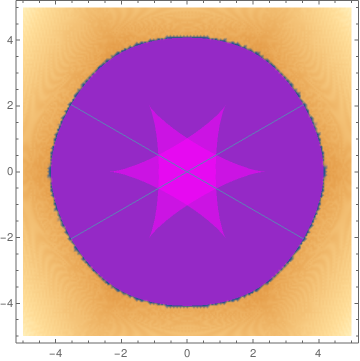}
        \caption{Phase locking for $\gamma_{ij}(0)=3$}
        \label{PhaseLock2}
    \end{minipage}%
\end{figure}

Figure (6) depicts the evolution for of the angles $\theta_1,\theta_2,\theta_3$ 
and coupling strengths $\gamma_{12},\gamma_{13},\gamma_{23}$ for a set of initial 
conditions that does not converge to a fixed point. In the numerical 
experiment in Figure (6) the parameter values are  $\alpha=.3$ and 
$\omega=(3,3,-6)/\sqrt{6}$ corresponding to $(a=0,b=1)$ in mean-zero frequency 
plane, with initial conditions $\theta_i(0)=0$ and $\gamma_{ij}(0)=1$. Refering to Figures (4) and (5) we can see that there is a 
stable fixed point for this frequency value but it is not reached from these
initial conditions. The graph of the left side of Figure (6) depicts the angles: it is clear that two oscillators synchronize with one another (the two with the same natural frequency), but they do not phase-lock to the third,
so the angular difference grows linearly with time. The lefthand graph in 
Figure (6) depicts  the coupling strengths as a function of time. One coupling strength (the one between the synchronized oscillators) quickly converges to a constant but the other two oscillate
in a periodic manner.

Figure (7) depicts the evolution for parameter values $\alpha=.3$ and $\omega=(-1,-.5,3/2)$ and  initial conditions $\theta_i(0)=0$ and $\gamma_{ij}(0)=1$, which lie in the basin of attraction for the phase=locked solution. In contrast to Figure (6) we see a brief transient followed by exponential convergence to a phase-locked solution. 
 
\begin{figure}[h]
\includegraphics[width=0.45\textwidth]{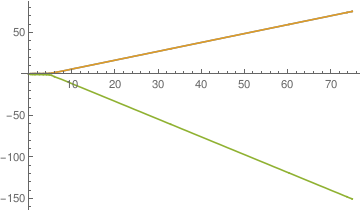} 
\includegraphics[width=0.45\textwidth]{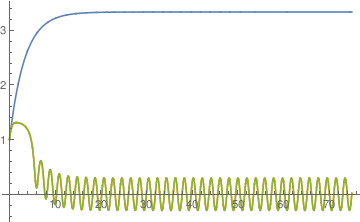}
\caption{ The angles $\theta_i(t)$ (left) and the edge weights $\gamma_{ij}(t)$ 
(right) for the three oscillator model in a non-phase-locked situation.} 
\end{figure}

\begin{figure}[h]
\includegraphics[width=0.45\textwidth]{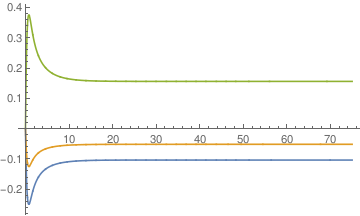} 
\includegraphics[width=0.45\textwidth]{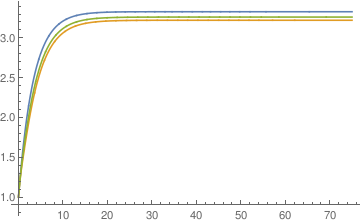}
\caption{ The angles $\theta_i(t)$ (left) and the edge weights $\gamma_{ij}(t)$ 
(right) for the three oscillator model in a phase-locked situation.} 
\end{figure}

\section{Conclusions}

We have shown that the fixed points of a Kuramoto model with Hebbian interactions are equivalent, under a simple map, to those of the classical Kuramoto model 
with a fixed interaction strength and, moreover, that the dimensions of the 
stable and unstable manifolds to these solutions are also the same. In the case of general (odd, periodic) interaction potential this is no longer the
case, but we are able to drastically reduce the dimensionality of the 
linearized flow by the use of a Schur complement identity.

{\bf Acknowledgements:} This was an undergraduate research project carried out under the auspices of the Illinois Geometry Lab (IGL). All of the authors would like to thank the IGL for support during the writing of this paper. J.C.B would additionally like to thank the National Science Foundation for support under 
award DMS 1211364, and Lee DeVille for many useful conversations.

\bibliography{Kuramoto}

\end{document}